\documentclass[11pt]{article}
\usepackage{epsfig}
\usepackage{amssymb}
\usepackage{amsmath}
\usepackage{amsfonts}

\usepackage[margin=1in]{geometry}
\usepackage{titlesec}
\usepackage[affil-it]{authblk}
\usepackage{breqn}
\usepackage{graphicx}

\usepackage{subfigmatrix}
\usepackage{placeins}
\usepackage{hyperref}
\usepackage{caption}

\newtheorem{theorem}{Theorem}

\newtheorem{corollary}{Corollary}

\usepackage{mathrsfs}

\newcommand{\Po}{\mathcal{P}}

\newcommand{\Pol}{\mathscr{P}}

\allowdisplaybreaks[2]

\begin{document}






%

\title{Experimental validation of volume-based comparison\\ for double-McCormick relaxations\footnote{Supported in part by ONR grant N00014-14-1-0315.}}
\date{\today}
\author{Emily Speakman}
\author{Han Yu}
\author{Jon Lee}

\affil{\rm Department of Industrial and Operations Engineering\\
\rm University of Michigan\\
\rm Ann Arbor, MI, U.S.A. \\ \quad\\
\rm $\{$eespeakm,~yuha,~jonxlee$\}$@umich.edu
}

\maketitle
\begin{abstract}
\noindent Volume is a natural geometric measure for comparing polyhedral relaxations of non-convex sets. Speakman and Lee gave volume formulae for comparing relaxations of trilinear monomials, quantifying the strength of various natural relaxations. Their work was motivated by the spatial branch-and-bound algorithm for factorable mathematical-programming formulations. They mathematically analyzed an important choice that needs to be made whenever three or more terms are multiplied in a formulation. We experimentally substantiate the relevance of their main results to the practice of global optimization, by applying it to difficult box cubic problems (boxcup). In doing so, we find that, using their volume formulae, we can accurately predict the quality of a relaxation for boxcups based on the (box) parameters defining the feasible region.

\end{abstract}

%
%


%
%

%
%


%
%

\thispagestyle{empty}
\clearpage
\setcounter{page}{1}

\section{Introduction}

\subsection{Measuring relaxations via volume in mathematical optimization}
In the context of mathematical optimization, there is often a natural tradeoff in the tightness
of a  \emph{convexification} (i.e., a convex relaxation) and the difficulty of optimizing over it. This idea was emphasized by
\cite{Lee_2007} in the context of mixed-integer nonlinear programming (MINLP) (see also the recent work \cite{Dey2015}). Of course this is also a well-known phenomenon for difficult 0/1
linear-optimization problems, where very tight relaxations are available via extremely heavy semidefinite-programming
relaxations (e.g., the Lasserre hierarchy), and the most effective relaxation for branch-and-bound/cut may well not be the tightest.
Earlier, again in the context of mathematical optimization,
\cite{Lee94} introduced the idea of using volume as a measure of
the tightness of a convex relaxation (for fixed-charge and vertex packing problem). Most of that mathematical work was asymptotic, seeking to
understand the quality of families of relaxations with a growing number of variables, but some of it was also substantiated experimentally in \cite{Lee_2007}.

\subsection{Spatial branch-and-bound}
The workhorse algorithm for global optimization of so-called factorable MINLPs (see \cite{McCormick76}) is spatial branch-and-bound (sBB) (see \cite{Adjiman98}, \cite{Ryoo96}, \cite{Smith99}).
sBB  decomposes model functions against a library of basic functions (e.g., $x_1x_2$, $x_1x_2x_3$, $x_1^{x_2}$, $\log(x)$, $\sin(x)$, $x^3$, $\sqrt{x}$, $\arctan(x)$).
We assume that each basic function
is a function of no more than 3 variables, and that we have a convex outer-approximation of the graph of each such basic function on
box domains. sBB realizes convexifications of  model functions by composing convexifications of basic functions.
sBB subdivides the domains of variables, and re-convexifies, obtaining stronger bounds on model functions.

\subsection{Using volume to guide decompositions for spatial branch-and-bound}
\cite{SpeakmanLee2015} 
applied the idea of \cite{Lee94}, but now in the context of the low-dimensional
relaxations of basic functions that arise in sBB. Specifically, \cite{SpeakmanLee2015} considered the basic function $f=x_1x_2x_3$ on box domains;
that is the graph
\[
\Po_h=\left\{ (f,x_1,x_2,x_3)\in\mathbb{R}^4 ~:~ f=x_1x_2x_3,~ x_i\in[a_i,b_i], i=1,2,3 \right\},
\]
 where $0\leq a_i <b_i$ are given constants. It is important to realize that $\Po_h$
 is relevant for any model where any three quantities (which could be complicated functions themselves) are multiplied. Furthermore, the case of nonzero lower bounds (i.e., $a_i>0$) is particularly relevant, especially when the multiplied
 quantities are complicated functions of model variables.
  Though polyhedral, $\Po_h$ has a rather complicated inequality description.
 Often lighter relaxations are used by modelers and MINLP software. \cite{SpeakmanLee2015} considered 3 natural relaxations of
 $\Po_h$. Thinking of the product as $x_1(x_2x_3)$ or $x_2(x_1x_3)$ or $x_3(x_1x_2)$, and employing the
 so-called McCormick relaxation twice, leads to three \emph{different} relaxations. \cite{SpeakmanLee2015}
 derive  analytic expressions for the volume of $\Po_h$ as well as all three of the natural relaxations.
 The expressions are formulae in the six constants $0\leq a_i <b_i$, $i=1,2,3$.
 In doing so, they quantify the quality of the various relaxations and provide recommendations for which to use.

\subsection{Our contribution}
The results of \cite{SpeakmanLee2015} are theoretical. Their utility for guiding modelers and
sBB implementers depends on the belief that volume is a good measure of
the quality of a relaxation. Morally, this belief is based on the idea that with no prior information on the form of an objective function, the solution of a relaxation
should be assumed to occur with a uniform density on the feasible region.
Our contribution is to experimentally validate the robustness of this theory in the context of a particular use case,
optimizing multilinear cubics over boxes (`boxcup'). There is considerable literature on techniques for optimizing quadratics, much of
which is developed and validated in the context of so-called `boxqp' problems, where we minimize $\sum_{i,j} q_{ij}x_i x_j$
over a box domain in $\mathbb{R}^n$. So our boxcup problems, for which we minimize $\sum_{i,j} q_{ijk}x_i x_j x_k$ over a box domain in $\mathbb{R}^n$, are   natural.

\cite{SpeakmanLee2015} find an ordering of the three natural relaxations by volume. But their formulae are for $n=3$.
Our experiments validate their theory as applied to our use case. We demonstrate that in the setting of `boxcup' problems, the
average objective discrepancy between relaxations very closely follows the prediction of the theory of \cite{SpeakmanLee2015}.
Moreover, we are able to demonstrate that these results are robust against sparsity of the cubic forms.

\cite{Lee94} defined the \emph{idealized radius} of a polytope in $\mathbb{R}^d$ as
essentially the $d$-th root of its volume  (up to some constants depending on $d$).
For a polytope that is very much like a ball in shape, we can expect that this
quantity is (proportional to) the ``average width'' of the polytope. The
average width arises by looking at `max minus min', averaged over all normalized linear objectives.
So, the implicit prediction of \cite{SpeakmanLee2015} is that the \emph{idealized radius}
should (linearly) predict the expected `max minus min' for normalized linear objectives.
We have validated this experimentally, and looked further into the \emph{idealized radial distance}
between pairs of relaxations, finding an even higher degree of linear association.

Finally, in the important case $a_1=a_2=0$, $b_1=b_2=1$, \cite{SpeakmanLee2015} found that
the two worst relaxations have the same volume, and the
greatest difference in volume between $\Po_h$ and the (two) worst relaxations occurs
when $a_3=b_3/3$. We present results of  experiments that clearly show that these
predictions via volume are again borne out on `boxcup' problems.

All in all, we present convincing experimental evidence that volume is a good predictor for
quality of relaxation in the context of sBB. Our results strongly suggest that the
theoretical results of \cite{SpeakmanLee2015} are important in devising decompositions
of complex functions in the context of factorable formulations and therefore our results help
inform both modelers and implementers of sBB.

\subsection{Literature review}

Computing the volume of a polytope is well known to be strongly \#P-hard
(see \cite{BW1991}). But in fixed dimension, or, in celebrated work,
by seeking an approximation via a randomized algorithm (see \cite{DFK1991}), positive results are available. Our work though is motivated not by algorithms for volume
calculation, but rather in certain situations where analytic formulae are available.

Besides  \cite{Lee94} and \cite{SpeakmanLee2015} (and the follow-on \cite{SpeakmanLee2016}), there have been a few papers on analytic formulae for volumes of polytopes that naturally arise in mathematical optimization; see \cite{KLS1997}, \cite{Steingr}, \cite{Burggraf2013}, \cite{Ardila2010}, \cite{Stanley1986}. But none of these works has attempted to apply their ideas to the low-dimensional polytopes that naturally arise in sBB, or even to apply their ideas to compare relaxations. One notable exception is \cite{CafieriLeeLiberti10},
which is a mostly computational precursor to \cite{SpeakmanLee2015}, focusing on
quadrilinear functions (i.e., $f=x_1x_2x_3x_4$).

There are many implementations of sBB. E.g., Baron \cite{BaronManual}, Couenne \cite{Belotti09}, Scip \cite{VigerskeGleixner2016} and Antigone \cite{misener-floudas:ANTIGONE:2014}.
Both Baron and Antigone use the complete linear-inequality description of $\Po_h$,
while Couenne and SCIP use an \emph{arbitrary} double-McCormick relaxation. Our results indicate that there are situations where the choice of Baron and Antigone is too heavy, and certainly even restricting to double-McCormick relaxations, Couenne and SCIP do not choose the best one.

There is a large literature on convexification of multilinear functions.
Most relevant to our work are:
the polyhedral nature of the convexification of the
graphs of  multilinear functions on box domains (see \cite{Rikun97});
the McCormick inequalities describing giving the complete linear-inequality description for bilinear functions on a box domain (see \cite{McCormick76});
 the complete linear-inequality description of $\Po_h$ (see \cite{Meyer04a} and \cite{Meyer04b}).

\section{Preliminaries}


\subsection{Convexifications}

Without loss of generality, we can relabel the three variables so that
\[\label{Omega}
a_1b_2b_3+b_1a_2a_3 \leq b_1a_2b_3 + a_1b_2a_3 \leq b_1b_2a_3 + a_1a_2b_3. \tag{$\Omega$}
\]
is satisfied by the variable bounds.  Given  $f:=x_1x_2x_3$ satisfying (\ref{Omega}),
there are three  choices of double-McCormick convexifications depending on the bilinear sub-monomial we convexify first. We could first group $x_1$ and $x_2$ and convexify $w=x_1x_2$; after this, we are left with the monomial $f=wx_3$ which we again convexify using McCormick. Alternatively, we could first  group variables $x_1$ and $x_3$, or variables $x_2$ and $x_3$.

To see how to perform these convexifications in general, we exhibit the double-McCormick convexification that first groups the variables $x_i$ and $x_j$.  Therefore we have $f=x_ix_jx_k$, and we let $w_{ij}=x_ix_j$, so $f=w_{ij}x_k$.  Convexifying $w_{ij}=x_ix_j$ using the standard McCormick inequalities, and then
convexifying $f=w_{ij}x_k$, again using the standard McCormick inequalities, we obtain the 8 inequalities:


\vspace{-7mm}
\begin{align*}
w_{ij}-a_jx_i-a_ix_j+a_ia_j \geq 0, &&f- a_kw_{ij} - a_ia_jx_k + a_ia_ja_k \geq 0,\\
-w_{ij} + b_jx_i + a_ix_j - a_ib_j \geq 0, &&-f +b_kw_{ij} + a_ia_jx_k - a_ia_jb_k \geq 0,\\
-w_{ij} +a_jx_i + b_ix_j - b_ia_j \geq 0, &&-f + a_kw_{ij} + b_ib_jx_k - b_ib_ja_k \geq 0,\\
w_{ij} - b_jx_i - b_ix_j + b_ib_j \geq 0, &&f - b_kw_{ij} - b_ib_jx_k + b_ib_jb_k \geq 0.
\end{align*}
\vspace{-3mm}

Using Fourier-Motzkin elimination (i.e., projection), we then eliminate the variable $w_{ij}$ to obtain the following system in our original variables $f, x_i,x_j$ and $x_k$.
\begin{align}
&x_i - a_i &\geq 0, \label{in1}\\
&x_j - a_j &\geq 0, \label{in2}\\
&f - a_ja_kx_i - a_ia_kx_j - a_ia_jx_k + 2a_ia_ja_k &\geq 0, \label{in3} \\
&f - a_jb_kx_i - a_ib_kx_j - b_ib_jx_k + a_ia_jb_k + b_ib_jb_k &\geq 0,\label{in4}\\
&- x_j + b_j &\geq 0, \label{in5}\\
&- x_i + b_i &\geq 0, \label{in6}\\
&f - b_ja_kx_i - b_ia_kx_j - a_ia_jx_k + a_ia_ja_k + b_ib_ja_k &\geq 0, \label{in7}\\
&f - b_jb_kx_i - b_ib_kx_j - b_ib_jx_k + 2b_ib_jb_k &\geq 0, \label{in8}\\
&-f + b_jb_kx_i + a_ib_kx_j + a_ia_jx_k - a_ia_jb_k - a_ib_jb_k &\geq 0, \label{in9}\\
&-f + a_jb_kx_i + b_ib_kx_j + a_ia_jx_k - a_ia_jb_k - b_ia_jb_k &\geq 0, \label{in10}\\
&- x_k + b_k &\geq 0, \label{in11}\\
&-f + b_ja_kx_i + a_ia_kx_j + b_ib_jx_k - a_ib_ja_k - b_ib_ja_k &\geq 0, \label{in12}\\
&-f + a_ja_kx_i + b_ia_kx_j + b_ib_jx_k - b_ia_ja_k - b_ib_ja_k &\geq 0, \label{in13}\\
&x_k - a_k &\geq 0, \label{in14} \\
&f - a_ia_jx_k & \geq 0, \label{in15} \\
&-f + b_ib_jx_k & \geq 0. \label{in16}
\end{align}

It is easy to see that the inequalities (\ref{in15}) and (\ref{in16}) are redundant:  (\ref{in15}) is $a_ja_k(\ref{in1}) + a_ia_k(\ref{in2}) + (\ref{in3})$, and (\ref{in16}) is $b_ja_k(\ref{in6}) + a_ia_k(\ref{in5}) +(\ref{in12})$.

We use the following notation in what follows.  For $i=1,2,3$, \emph{system} $i$
is defined to be the system of inequalities obtained by first grouping the
pair of variables $x_j$ and $x_k$, with $j$ and $k$ different from $i$.
$\Po_i$ is defined to be the solution set of this system.

\subsection{Volumes of convexifications}

For a convex body $C \subset\mathbb{R}^d$, we denote its \emph{volume} (i.e.,  Lebesgue measure) by $\text{vol}(C)$.
The main results of \cite{SpeakmanLee2015} are as follows.
\begin{theorem}[\cite{SpeakmanLee2015}]
\label{TheoremPH}
\begin{dmath*}
\emph{vol}(\Po_h) = (b_1-a_1)(b_2-a_2)(b_3-a_3)\times\\
\phantom{a}\qquad\qquad\left(b_1(5b_2b_3 - a_2b_3 - b_2a_3 - 3a_2a_3) + a_1(5a_2a_3 - b_2a_3 - a_2b_3 -3b_2b_3)\right)/24.
\end{dmath*}
\end{theorem}

\begin{theorem}[\cite{SpeakmanLee2015}]
\label{TheoremP1}

\begin{dmath*}
\emph{vol}(\Po_1) = \emph{vol}(\Po_h) + (b_1-a_1)(b_2-a_2)^2(b_3-a_3)^2\times\\
\phantom{a}\qquad\qquad\qquad\qquad\frac{3(b_1b_2a_3 - a_1b_2a_3 + b_1a_2b_3 -a_1a_2b_3) + 2(a_1b_2b_3 - b_1a_2a_3)}{24(b_2b_3-a_2a_3)}.
\end{dmath*}
\end{theorem}

\begin{theorem}[\cite{SpeakmanLee2015}]
\label{TheoremP2}

\begin{dmath*}
\emph{vol}(\Po_2) = \emph{vol}(\Po_h) +
\frac{(b_1-a_1)(b_2-a_2)^2(b_3-a_3)^2\left(5(a_1b_1b_3-a_1b_1a_3) + 3(b_1^2a_3 - a_1^2b_3)\right)}{24(b_1b_3-a_1a_3)}.
\end{dmath*}
\end{theorem}

\begin{theorem}[\cite{SpeakmanLee2015}]
\label{TheoremP3}
\begin{dmath*}
\emph{vol}(\Po_3) = \emph{vol}(\Po_h)+
\frac{(b_1-a_1)(b_2-a_2)^2(b_3-a_3)^2\left(5(a_1b_1b_2-a_1b_1a_2) + 3(b_1^2a_2 -a_1^2b_2)\right)}{24(b_1b_2-a_1a_2)}.
\end{dmath*}
\end{theorem}

\begin{corollary}[\cite{SpeakmanLee2015}]
\label{TheoremCompare}
\[
\emph{vol}(\Po_h) \leq \emph{vol}(\Po_3) \leq \emph{vol}(\Po_2) \leq \emph{vol}(\Po_1).
\]
\end{corollary}

\begin{corollary}[\cite{SpeakmanLee2015}]
\label{TheoremWorst}
For the special case of $a_1=a_2=0$, $b_1=b_2=1$,  and fixed $b_3$, the greatest difference in volume for $\Po_3(=\Po_h$) and $\Po_2$ (or $\Po_1$)  occurs when $a_3=b_3/3$.
\end{corollary}

\subsection{From volumes to gaps}

Volume seems like an awkward measure to compare relaxations, when typically we are interested in objective-function gaps.
Following \cite{Lee94}, the \emph{idealized radius} of a convex body $C \subset\mathbb{R}^d$ is
\[
\rho(C):= \left(\mbox{vol}(C)/\mbox{vol}(B_d)\right)^{1/d},
\]
where $B_d$ is the (Euclidean) unit ball in $\mathbb{R}^d$. $\rho(C)$ is simply the radius of a ball having the same volume as $C$.
The \emph{idealized radial distance} between convex bodies $C_1$ and $C_2$ is simply $|\rho(C_1)-\rho(C_2)|$.
If $C_1$ and $C_2$ are concentric balls,  say with $C_1\subset C_2$, then the idealized radial distance
between them is the (radial) height of $C_2$ above $C_1$.
The \emph{mean semi-width} of $C$ is simply
\[
\frac{1}{2}\int_{\|c\|=1}
\left(
 \max_{x\in C} c'x - \min_{x\in C} c'x
 \right)
  d\psi,
\]
where $\psi$ is the $(d - 1)$-dimensional Lebesgue measure on the boundary of $B_d$, normalized
so that $\psi$ on the entire boundary is unity. If $C$ is itself a ball, then
(i) its idealized radius is in fact its radius, and (ii) its width in any unit-norm direction $c$ is constant,
and so (iii) its (idealized) radius is equal to its mean semi-width.

{\bf Key point:} What we can hope is that our relaxations are round enough so that choosing one of small volume (which is proportional to the idealized radius) is a good proxy for choosing the relaxation by \emph{mean width} (which is the same as mean objective-value range).

\section{Computational Experiments}

\subsection{Box cubic programs and 4 relaxations}

Our experiments are aimed at the following natural problem which is concerned with optimizing a linear function on trinomials. Let $H$ be a 3-uniform hyper-graph on $n$ vertices.
Each hyper-edge of $H$ is a set of 3 vertices, and we denote the set of
hyperedges by $E(H)$. If $H$ is complete, then $|E(H)|=\binom{n}{3}$.
We associate with each vertex $i$ a variable $x_i\in[a_i,b_i]$,
and with each hyper-edge $\{i,j,k\}$ the trinomial $x_ix_jx_k$
and a coefficient $q_{\tiny ijk}$ ($1\leq i<j<k\leq n$).
We now formulate the associated \emph{boxcup} (`box cubic problem'):
\vskip-10pt
\begin{equation}
\label{origprob}
 \min_{x\in\mathbb{R}^n}\left\{ \sum_{\{i,j,k\}\in E(H)} q_{\tiny ijk}\,   x_i\,   x_j\,  x_k ~:~
x_i \in[a_i,b_i],~ i=1,2,\ldots,n\right\}. \tag{BCUP}
\end{equation}
%
 The name is in analogy with the well-known \emph{boxqp}, where just two terms (rather than three) are multiplied  (`box' refers to the feasible region and `qp' refers to `quadratic program).

(\ref{origprob}) is a difficult nonconvex global-optimization problem.  Our goal here is not to solve instances of this problem, but rather to solve a number of different \emph{relaxations} of the problem and see how the results of these experiments correlate with the volume results of \cite{SpeakmanLee2015}. In this way, we seek to
determine if the guidance of \cite{SpeakmanLee2015} is relevant to
modelers and those implementing sBB.

We have seen how for a single trilinear term $f=x_ix_jx_k$, we can build four distinct relaxations: the convex hull of the feasible points, $\Po_h$, and three relaxations arising from double McCormick: $\Po_1$, $\Po_2$ and $\Po_3$.  To obtain a relaxation of (\ref{origprob}), we choose a relaxation $\Po_{\ell}$, for some $\ell=1,2,3,h$ and apply this same relaxation method to \emph{each} trinomial of (\ref{origprob}).  We therefore obtain 4 distinct linear relaxations of the form:
\vskip-10pt
\begin{equation*}
\min_{(x,f)\in \Pol_{\ell}}\left\{ \sum_{\{i,j,k\}\in E(H)} q_{\tiny ijk} \,  f_{\tiny ijk}  
 \right\}.
\end{equation*}
where $\Pol_{\ell}$, $\ell=1,2,3,h$ is the polytope in dimension ${n\choose 3}+ n$ arising from using relaxation $\Po_{\ell}$ on each trinomial.  This linear relaxation is a linear inequality system involving the $n$ variables $x_i$ ($i=1,2,\ldots,n$), and the ${n \choose 3}$ new `function variables' $f_{ijk}$. These `function variables' model $x_i\, x_j\, x_k$.

For our experiments, we randomly generate box bounds $[a_i,b_i]$ on  $x_i$, for each $i=1,\ldots,n$ independently, by choosing (uniformly) random pairs of integers $0\leq a_i < b_i\leq 10$. With each realization of these bounds, we get relaxation feasible regions $\Po_{\ell}$, for $\ell=1,2,3,h$.

\subsection{3 scenarios for $Q$}

We have 3 scenarios for the hypergraph $H$ of (\ref{origprob}), all with $|E(H)|=20$ monomials:
\begin{itemize}
\item Our {\bf dense} scenario has $H$ being a \emph{complete} 3-uniform hypergraph
on $n=6$ vertices ($\binom{6}{3}=20$). We note that each of the $n=6$ variables
appears in $\binom{6-1}{3-1}=10$ of the 20 monomials, so there is
considerable overlap in variables between trinomials.
\item Our {\bf sparse} scenario has hyperedges: $\{1,2,3\}$,
$\{2,3,4\}$, $\{3,4,5\}$ \ldots $\{18,19,20\}$, $\{19,20,1\}$, $\{20,1,2\}$.
Here we have $n=20$ variables and each variable is in only 3 of the trinomials.
\item Our {\bf very sparse} scenario has $n=30$ variables and each variable is in only 2 of the trinomials.
\end{itemize}
For each scenario, we generate 30  \emph{sets} of  bounds $[a_i,b_i]$ on  $x_i$ ($i=1,\ldots,n$).
To control the variation in our results, and considering that the
scaling of $Q$ is arbitrary, we generate 100,000 random  $Q$ with $|E(H)|$ entries, uniformly distributed on the unit
sphere in $\mathbb{R}^{|E(H)|}$. 
Then, for each $Q$, we both minimize and maximize $\sum_{i<j<k} q_{\tiny ijk} \,  f_{\tiny ijk}$ over each $\Pol_{\ell}, \ell=1,2,3,h$ and each set of bounds.


\subsection{Quality of relaxations}

For each $Q$ we take the difference in the optimal values, i.e. the maximum value minus the minimum value; this can be thought of as the width of the polytope in the direction $Q$.  We then average these widths for each $\Pol_{\ell}$, $\ell=1,2,3,h$, across the 100,000 realizations of $Q$ (which results in very small standard errors),
 and we refer to this quantity $\omega(\Pol_{\ell})$ as
 the \emph{quasi mean width} of the relaxation. It is not quite the geometric mean width, because we do not have objective terms for all variables in (\ref{origprob}) (i.e., we have no objective terms $\sum_{i=1}^n c_ix_i$).

We seek to investigate how well the volume formulae of \cite{SpeakmanLee2015}, comparing the volumes of the polytopes $\Po_{\ell}$ ($\ell=1,2,3,h$),
can be used to predict the quality of the relaxations $\Pol_{\ell}$ ($\ell=1,2,3,h$) as measured by their quasi mean width.

\begin{figure}[p]
    \begin{subfigmatrix}{1}
  \subfigure[dense case]{\includegraphics[width=0.6\textwidth,keepaspectratio]{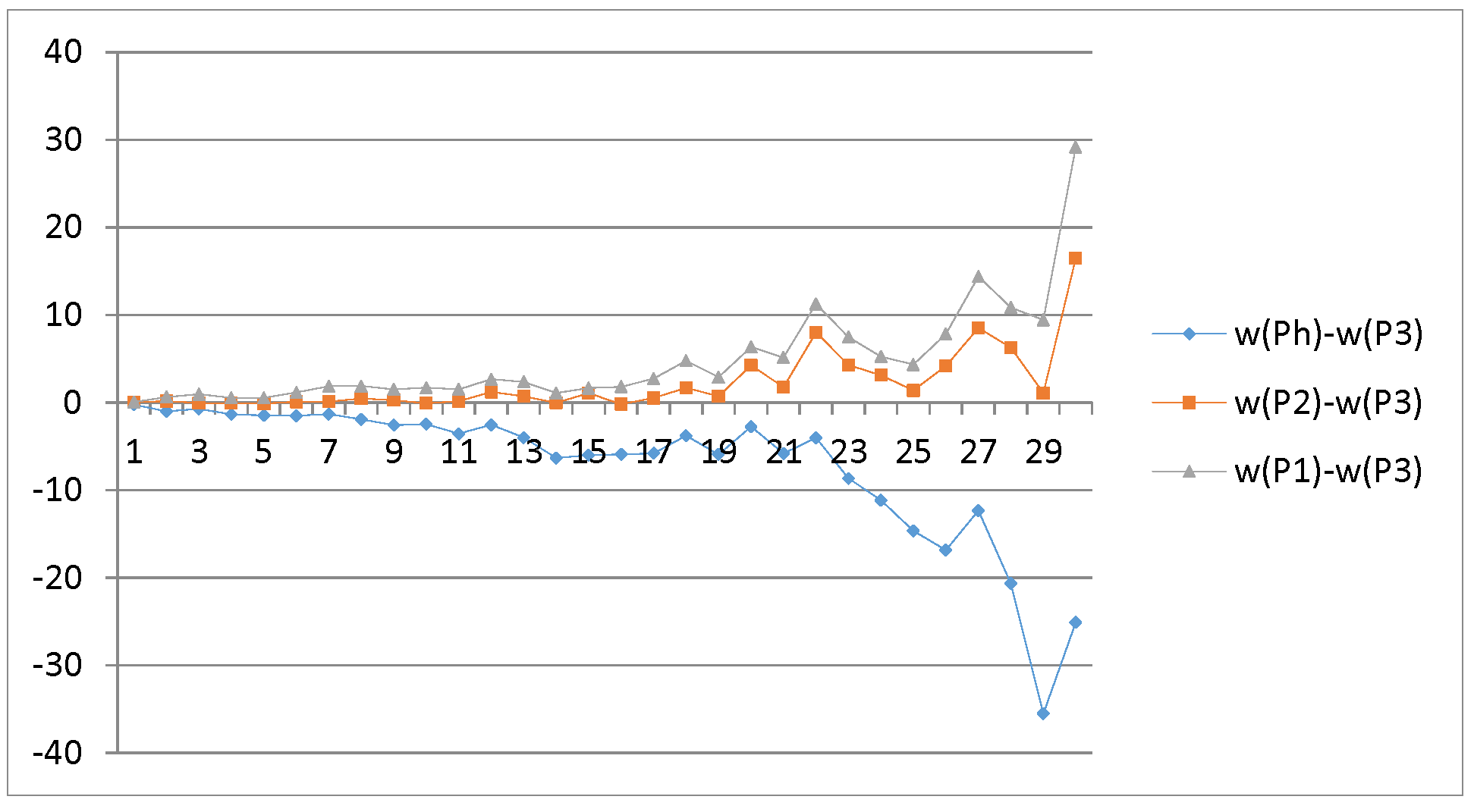}}
  \subfigure[sparse case]{\includegraphics[width=0.6\textwidth,keepaspectratio]{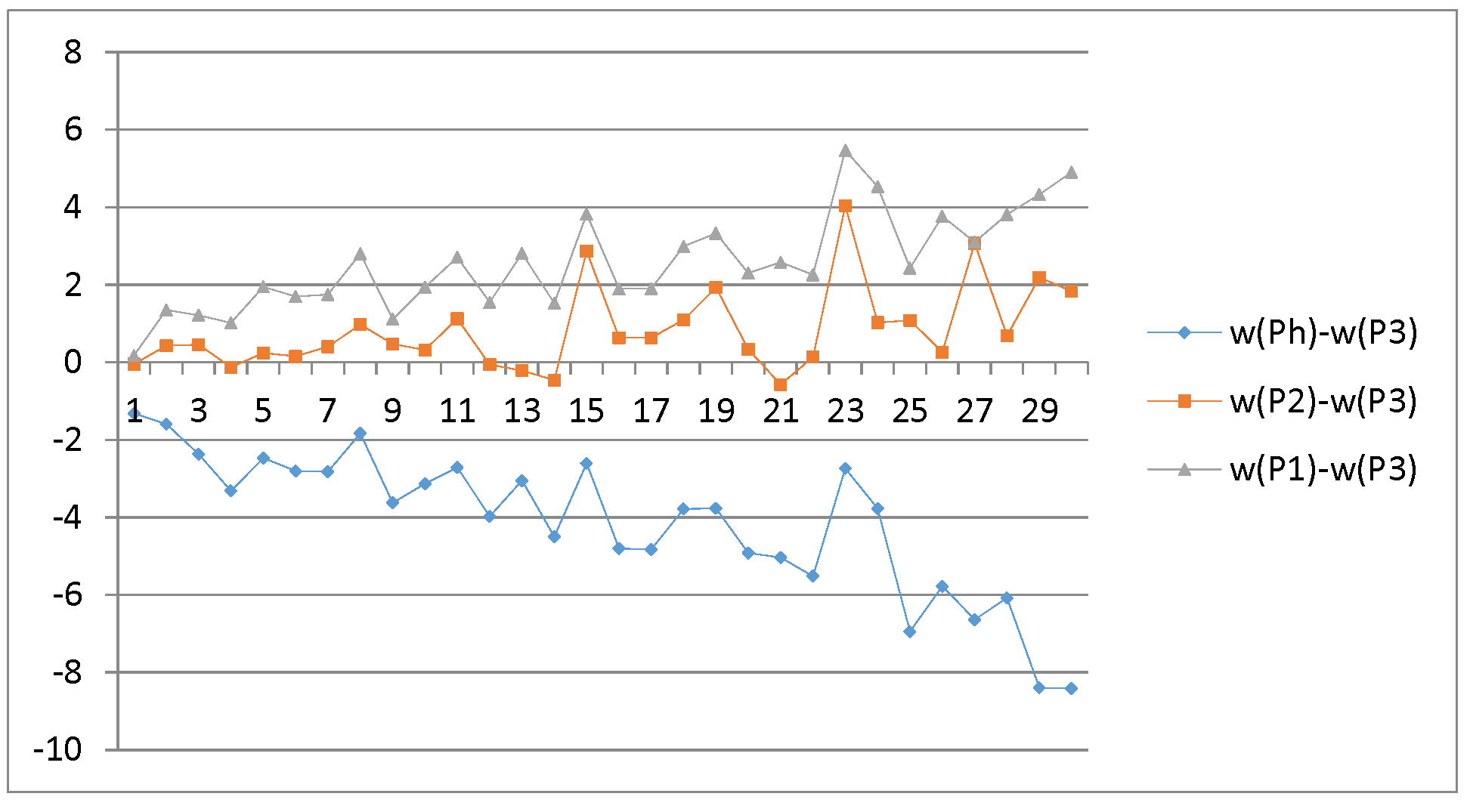}}
  \subfigure[very-sparse case]{\includegraphics[width=0.6\textwidth,keepaspectratio]{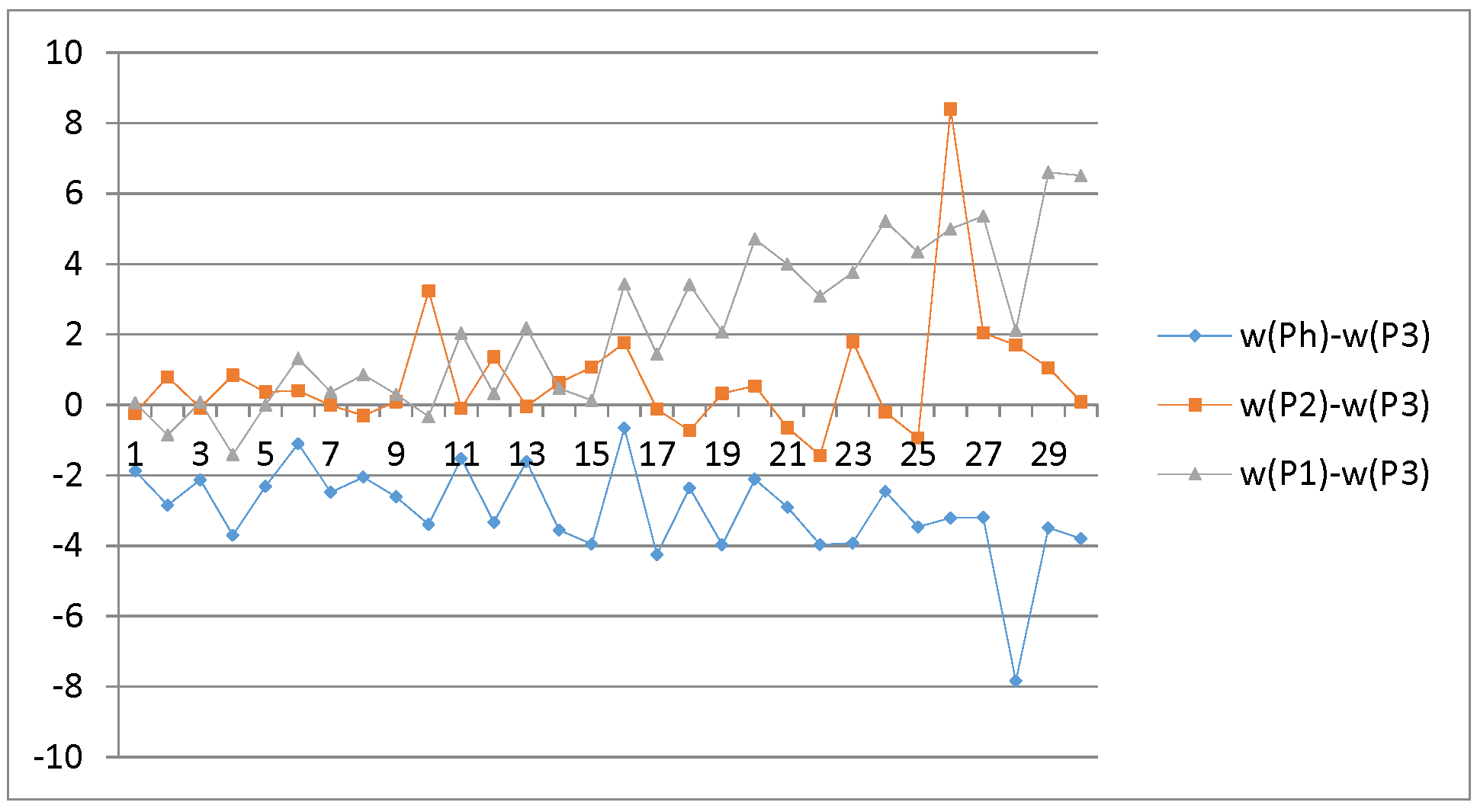}}
    \end{subfigmatrix}
    \caption{Quasi-mean-width differences}
    \label{fig:voldiffs}
 \end{figure}

Figure \ref{fig:voldiffs} consists of a plot for each scenario:  dense,  sparse, and  very sparse.  Each plot illustrates the difference in quasi mean width  between $\Pol_3$ (using the `best' double McCormick) and each of the other relaxations.  Each point represents a choice of bounds and the instances are sorted by $\omega(\Pol_1)-\omega(\Pol_h)$.  In all three plots $\omega(\Pol_h)-\omega(\Pol_3)$ is non-positive, which is to be expected because $\Pol_h$ is contained in each of the three double-McCormick relaxations.  Furthermore the plots illustrate that the general trend is for $\omega(\Pol_2)-\omega(\Pol_3)$ and $\omega(\Pol_1)-\omega(\Pol_3)$ to be positive and also for $\omega(\Pol_1)-\omega(\Pol_3)$ to be greater than $\omega(\Pol_2)-\omega(\Pol_3)$.  This agrees with Corollary \ref{TheoremCompare} and gives strong validation for the use of volume to measure the strength of different relaxations.  It confirms that given a choice of the double-McCormick relaxations, $\Pol_3$ is the one to choose.

However, there are a few exceptions to the general trend and these exceptions are most apparent in the very sparse case.  In both the dense case and the sparse case we only see a deviation from the trend on a small number of occasions when $\Pol_2$ is very slightly better than $\Pol_3$.  In each of these cases, the difference seems to be so small that we can really regard $\Pol_3$ and $\Pol_2$ as being equivalent from a practical viewpoint.  In the very sparse case, the general trend is still followed, but we see a few more cases where $\Pol_2$ is slightly better than $\Pol_3$.  We also see that in a few instances, $\Pol_1$ is better than $\Pol_2$ and occasionally even slightly better than $\Pol_3$.

However it is important to note that when we consider the sparse and very sparse cases, the differences in quasi mean width between \emph{any} two of the relaxations is much smaller than these differences in the dense case.   If we were to take the sparsity of $H$ to the extreme and run our experiments with $n=60$ and each variable only in one trinomial, the difference in quasi mean width between \emph{any} two of the polytopes will become zero for these boxcup problems.  Therefore, it is not surprising that our results diverge from the general trend as $H$ becomes sparser.

\begin{figure}[p]
    \begin{subfigmatrix}{1}
  \subfigure[dense case]{\includegraphics[width=0.5\textwidth,keepaspectratio]{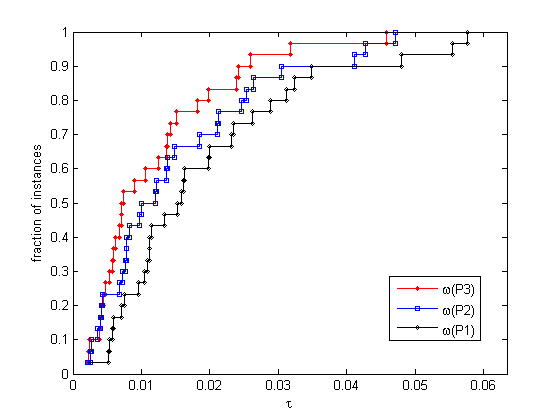}}
  \subfigure[sparse case]{\includegraphics[width=0.5\textwidth,keepaspectratio]{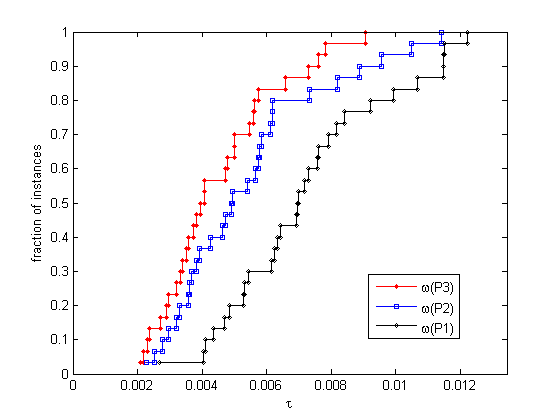}}
  \subfigure[very-sparse case]{\includegraphics[width=0.5\textwidth,keepaspectratio]{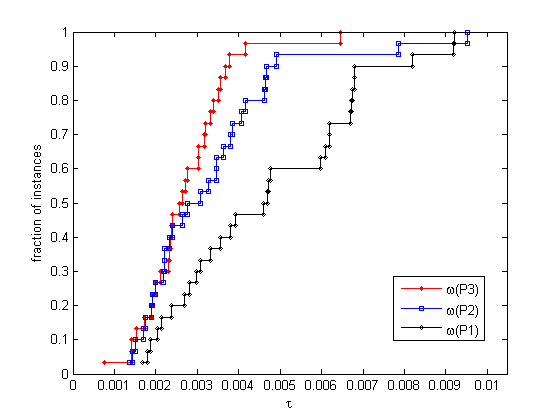}}
    \end{subfigmatrix}
    \caption{Quasi-mean-width performance profiles}
    \begin{center} Displays the fraction of instances where the quasi mean width of $\Pol_{\ell}$ is within\\ a factor $\alpha=e^\tau$ of the mean width of $\Pol_h$.
    Note that for small $\tau$, $e^\tau \approx 1+\tau$.\end{center}
    \label{fig:volpps}
 \end{figure}

Using the common technique of `performance profiles' (see \cite{Dolan2002}), we can illustrate the differences in quasi mean width of the three double-McCormick relaxations in another way.  We obtained the matlab code ``perf.m"  which was adapted to create these plots from the link contained in \cite{Sofi2015}.  
 Figure \ref{fig:volpps} shows a performance profile for each of the dense, sparse and very sparse scenarios.  For each choice of bounds, $\Pol_h$ gives the least quasi mean width (because it is contained in each of the other relaxations).  Our performance profiles display the fraction of instances where the quasi mean width of $\Pol_{\ell}$ is within a factor $\alpha$ of the mean width of $\Pol_h$, for $\ell =1,2,3$.  The plots are natural log plots where the horizonal axis is $\tau:=\ln(\alpha)$.  Using this measure, we see that the trend in \emph{all} cases is that $\Pol_3$ dominates $\Pol_2$ which in turn dominates $\Pol_1$.  In the very sparse case, we see that $\Pol_3$ and $\Pol_2$ are very close for small factors $\alpha$.  In general, all three relaxations are within a small factor of the hull.  Displaying the results in this manner gives us a way to see quickly which relaxation performs best for the majority of instances.  Again, we see agreement with the prediction of Corollary \ref{TheoremCompare} and confirmation that $\Pol_3$ is the best double-McCormick relaxation.

\subsection{Validating the relationship between volume and mean width}

\begin{figure}[p]
    \begin{subfigmatrix}{2}
  \subfigure[$\Pol_h$]{\includegraphics[width=0.47\textwidth,keepaspectratio]{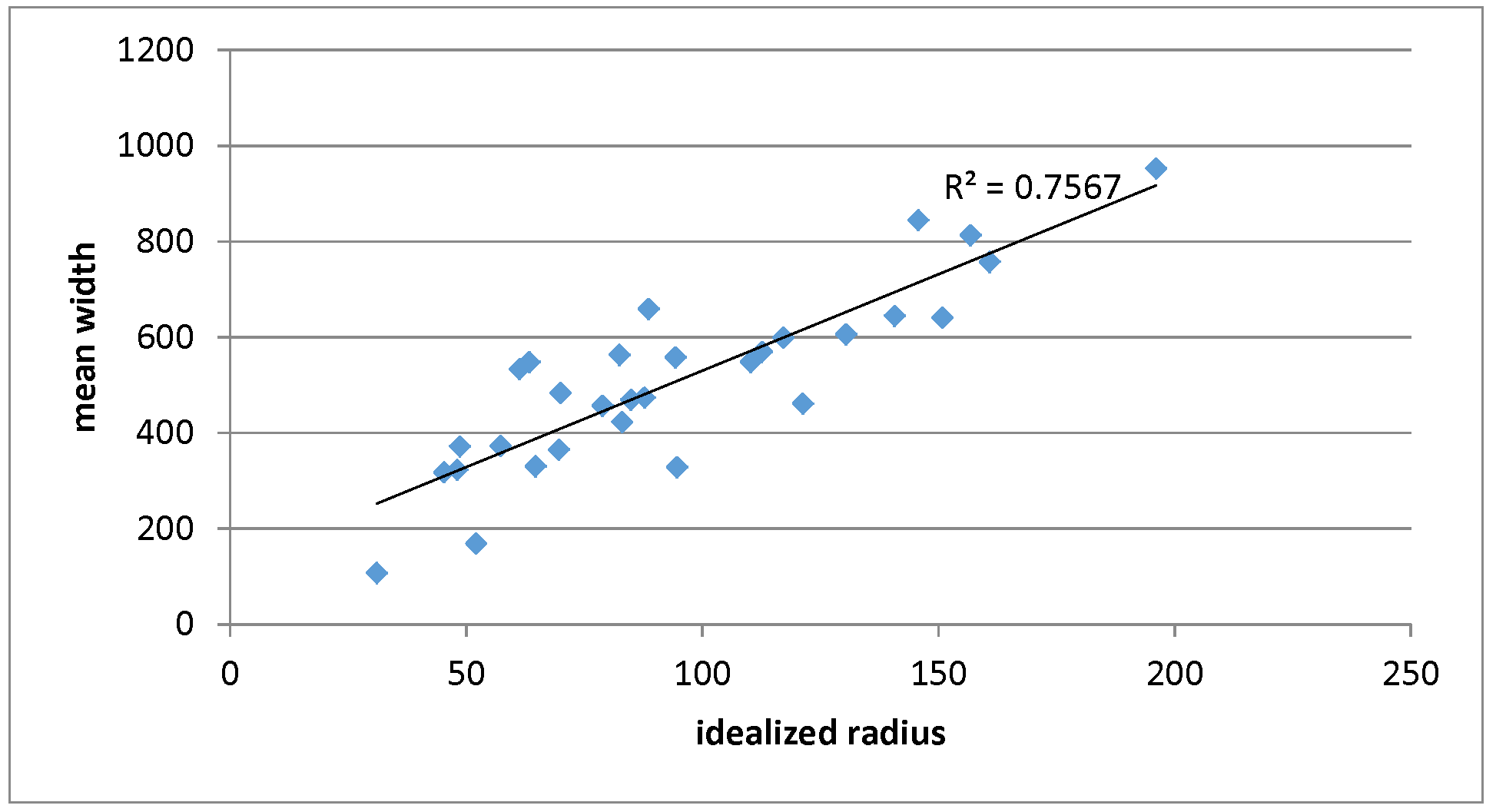}}
  \subfigure[$\Pol_3$]{\includegraphics[width=0.47\textwidth,keepaspectratio]{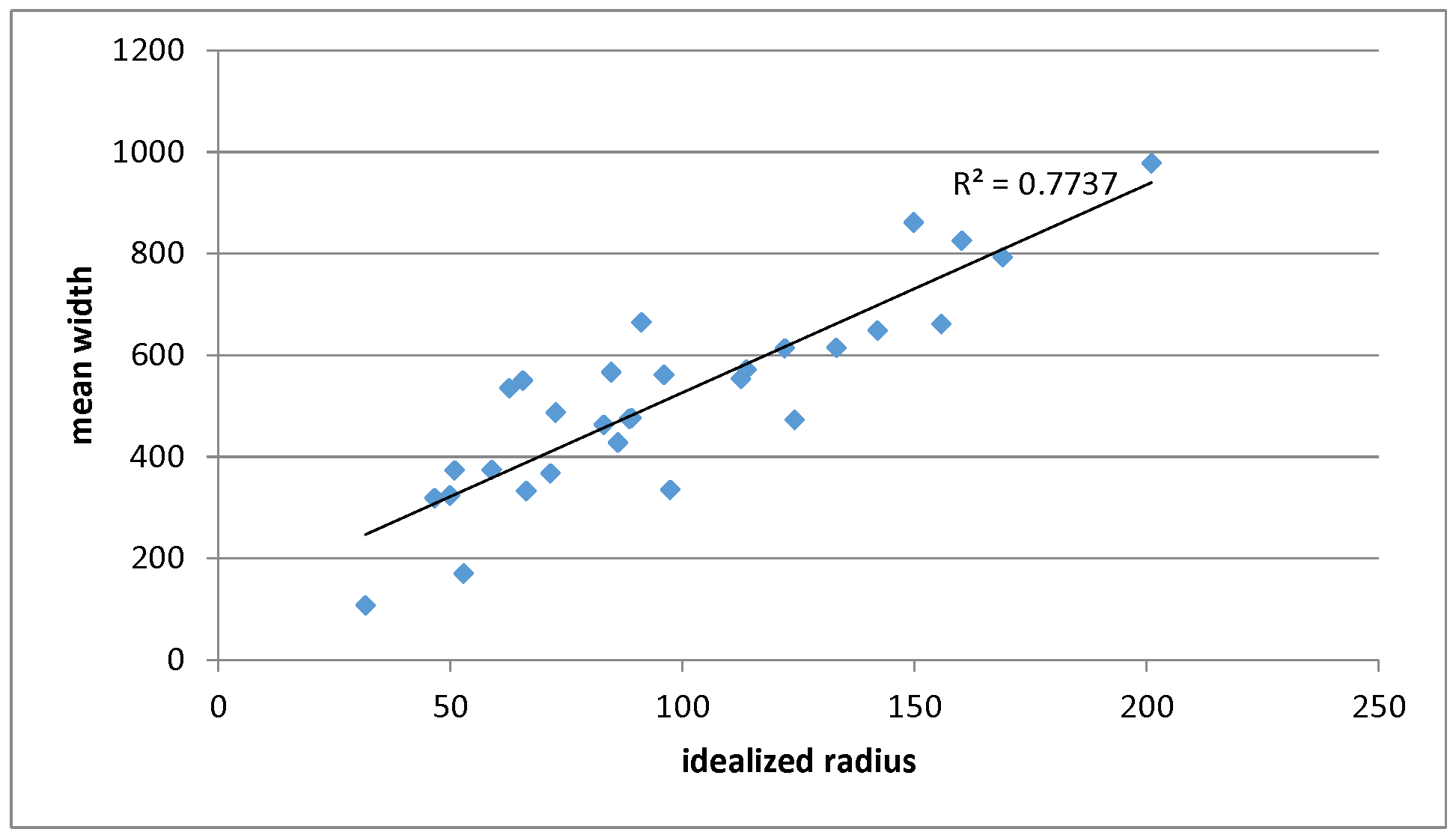}}
  \subfigure[$\Pol_2$]{\includegraphics[width=0.47\textwidth,keepaspectratio]{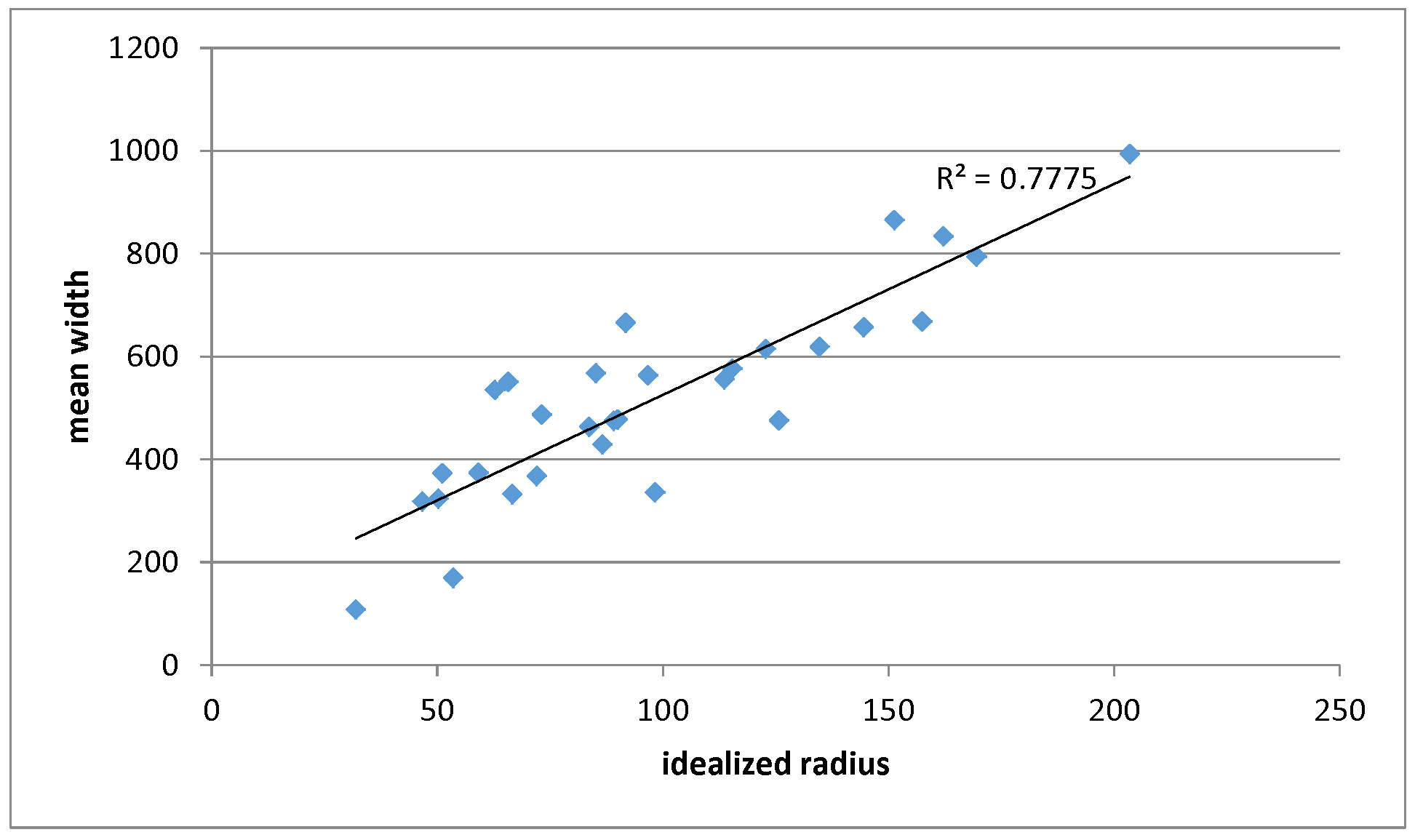}}
  \subfigure[$\Pol_1$]{\includegraphics[width=0.47\textwidth,keepaspectratio]{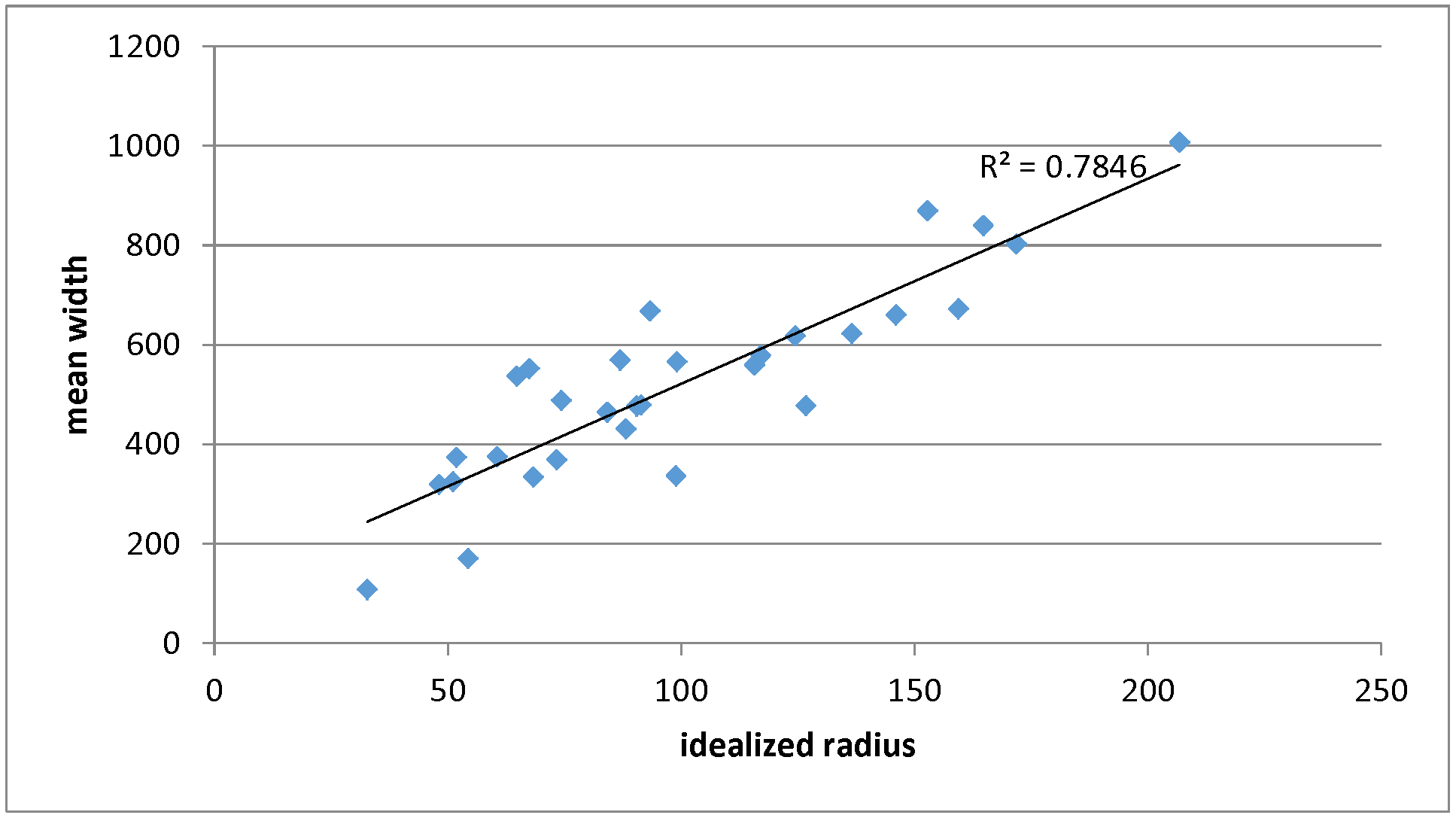}}
    \end{subfigmatrix}
    \caption{idealized radius predicting quasi mean width}
    \label{fig:widths}
 \end{figure}

Using the \cite{SpeakmanLee2015} formulae, we calculate the volume of the relaxation for each individual trinomial, $\Po_{\ell}$. We then we take the fourth root of these volumes and sum over all $|E(H)|$ trinomials to obtain a kind of `aggregated idealized radius' for each relaxation and each set of bounds. Restricting our attention to the dense scenario, in Figure \ref{fig:widths},
we compare these aggregated idealized radii with quasi mean width, across all relaxations
$\Pol_{\ell}, \ell=1,2,3,h$ and each set of bounds (each point in each scatter plot corresponds to a choice of bounds).
We see a high $R^2$ coefficient in all cases,
so we may conclude that volume really is a good predictor of relaxation width.

\begin{figure}[p]
    \begin{subfigmatrix}{2}
  \subfigure[$\Pol_3$ vs. $\Pol_h$]{\includegraphics[width=0.47\textwidth,keepaspectratio]{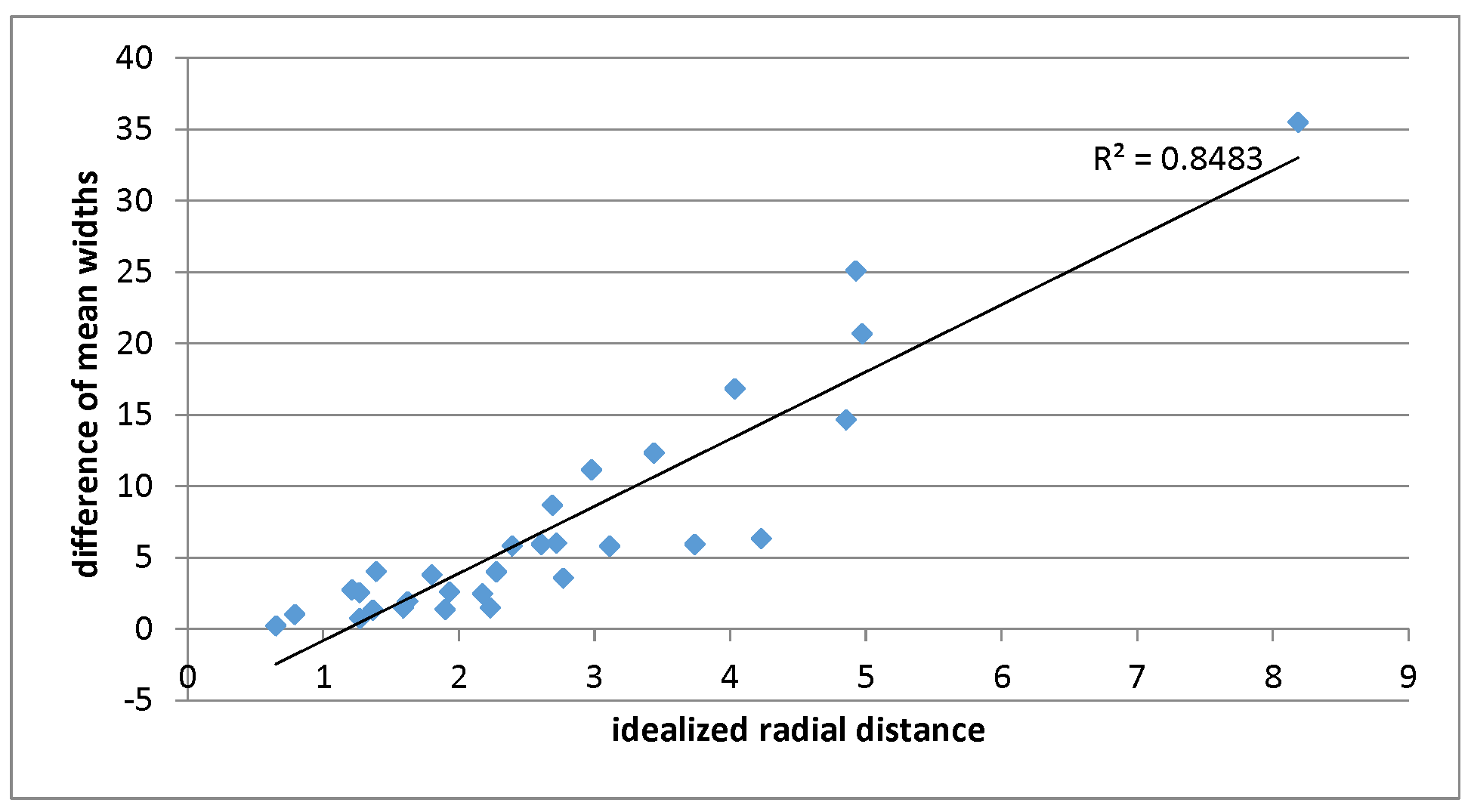}}
  \subfigure[$\Pol_2$ vs. $\Pol_h$]{\includegraphics[width=0.47\textwidth,keepaspectratio]{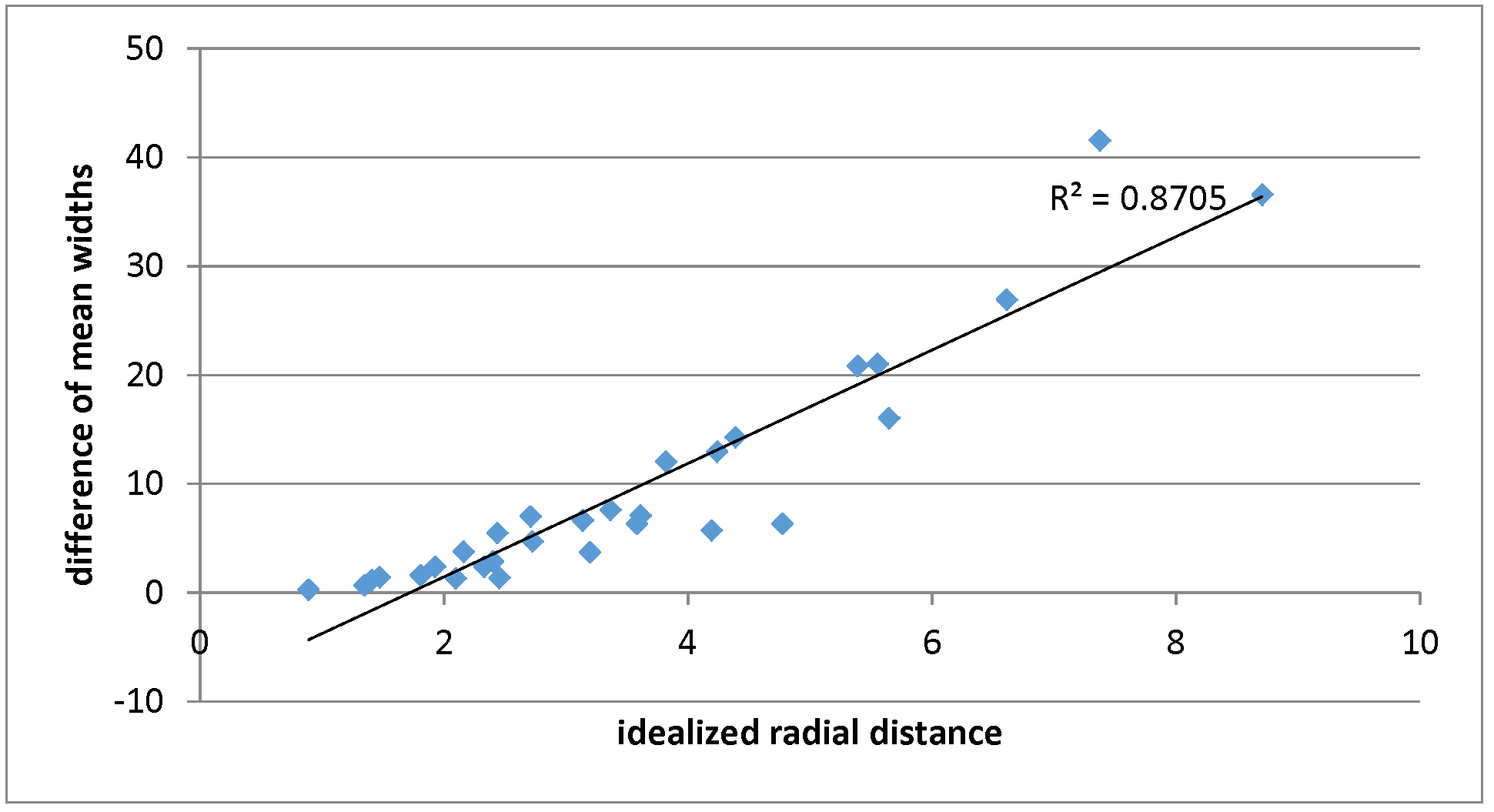}}
  \subfigure[$\Pol_1$ vs. $\Pol_h$]{\includegraphics[width=0.47\textwidth,keepaspectratio]{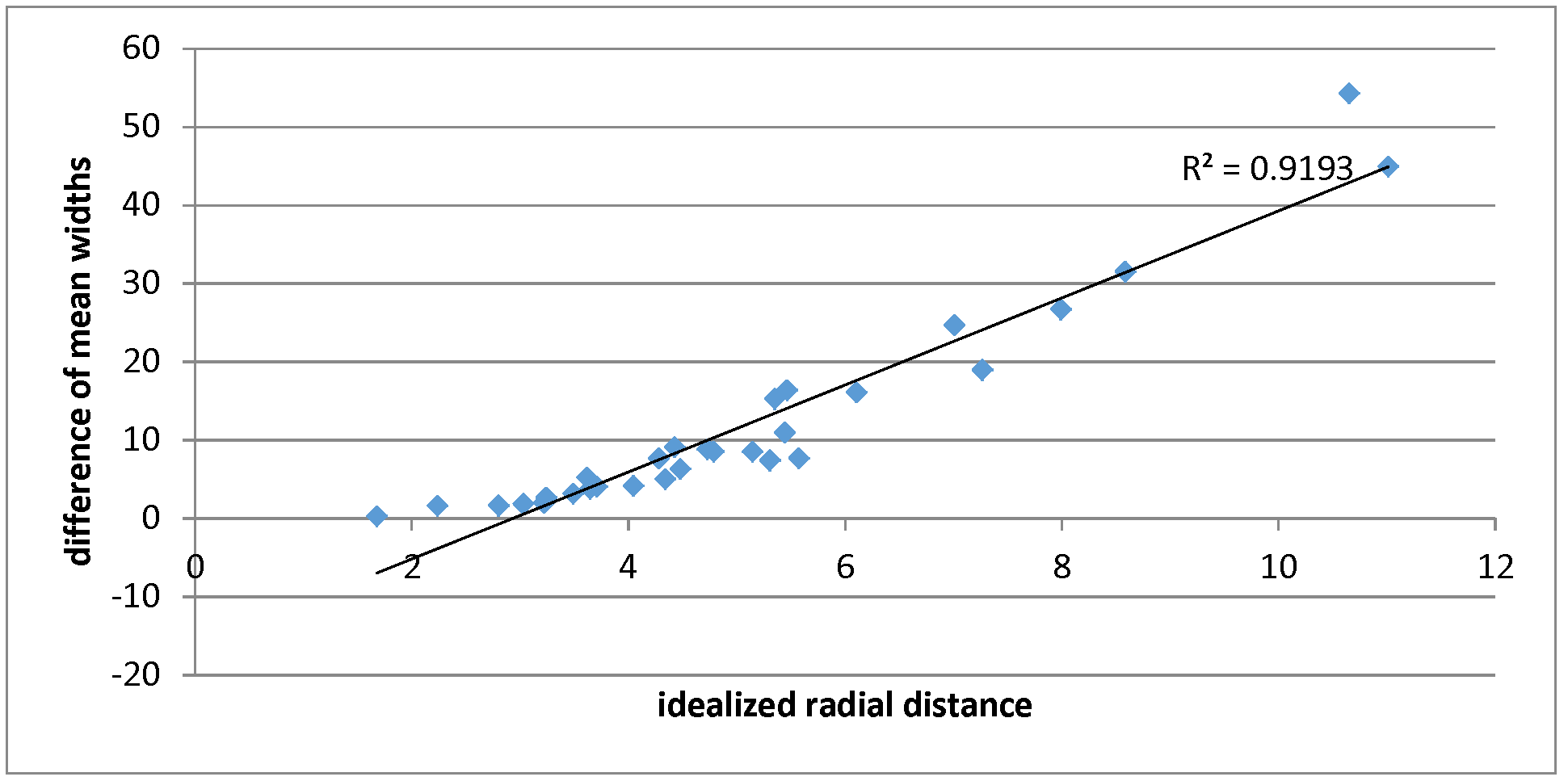}}
    \end{subfigmatrix}
    \caption{idealized radial distance predicting quasi mean width difference}
    \label{fig:widthdiffs}
 \end{figure}

We also compute the difference in width between polytope pairs: $\Pol_h$ and $\Pol_3$, $\Pol_3$ and $\Pol_2$, $\Pol_2$ and $\Pol_1$ for each direction $Q$.  We then average these width differences for each polytope and each set of bounds, across the 100,000 realizations of $Q$. We refer to this result as the \emph{quasi mean width difference} of the pair of polytopes. In Figure \ref{fig:widthdiffs}, we similarly
compare aggregated idealized radial differences with quasi mean width differences.
We see even higher $R^2$ coefficients, validating volume as an excellent predictor
of average objective gap between pairs of relaxations.

\subsection{A worst case}

Our final set of experiments relate to a `worst case' as described in \cite{SpeakmanLee2015}.  In the important special case of $a_1=a_2=0$ and $b_1=b_2=1$, the two `bad' double-McCormick relaxations have the same volume and the `good' double McCormick is exactly the hull.  In addition, the greatest difference in volume between the hull and the bad relaxations occurs when $a_3= \frac{b_3}{3}$.

We compute the same results as we have discussed before (i.e. the differences in quasi mean width between the relaxations) with $n=6$, but now instead of using random bounds, we fix $a_1=a_2=0$ and $b_1=b_2=1$.  We also fix $b_3$ and run the experiments for $a_3=1,2,\dots, b_3-1$.  Here, we only consider the ${5 \choose 3}=10$ trinomials that have the form $x_jx_kx_6$.

\begin{figure}[h!]
    \begin{subfigmatrix}{2}
  \subfigure[$b_3=30$]{\includegraphics[width=0.47\textwidth,keepaspectratio]{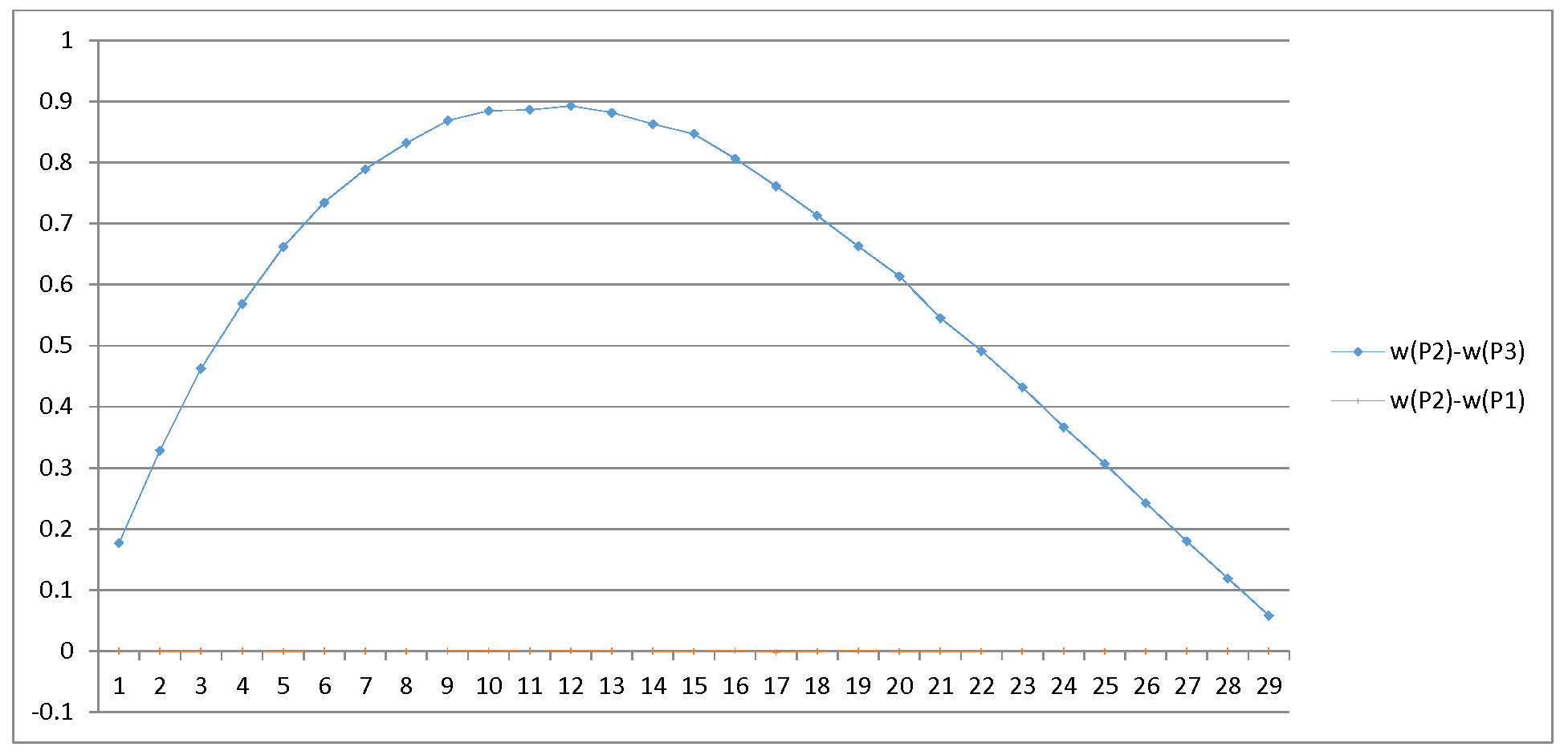}}
  \subfigure[$b_3=60$]{\includegraphics[width=0.47\textwidth,keepaspectratio]{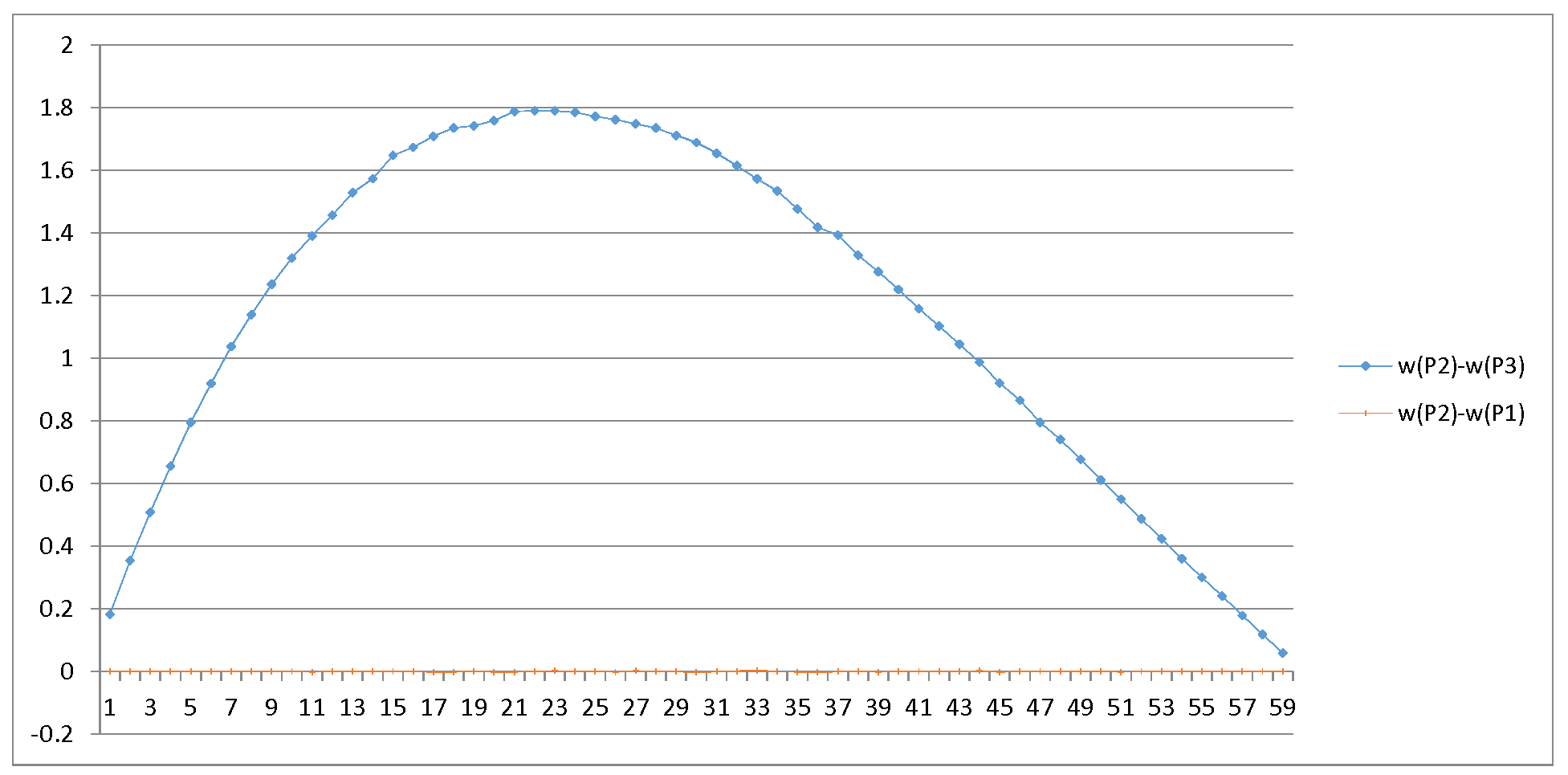}}
 \subfigure[$b_3=90$]{\includegraphics[width=0.47\textwidth,keepaspectratio]{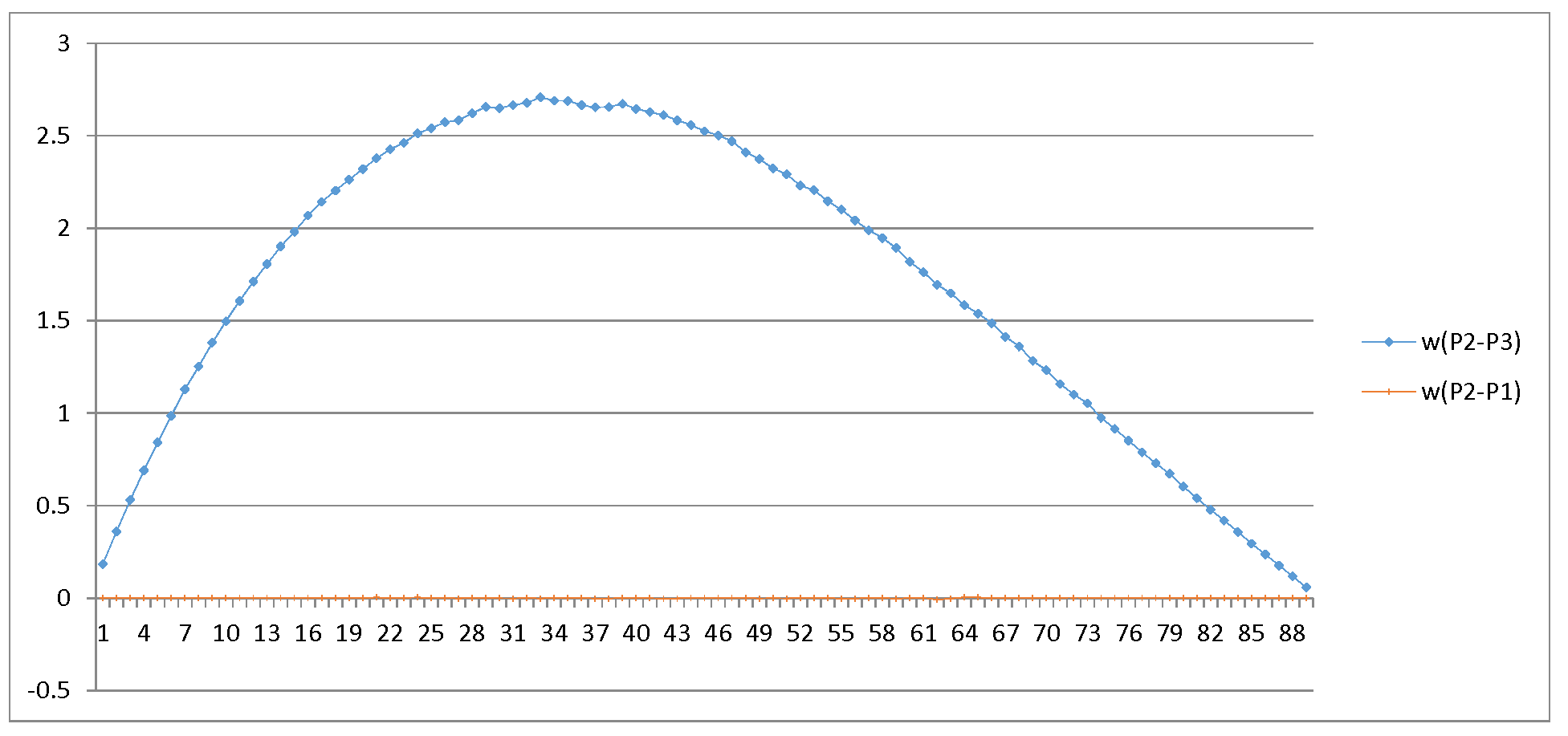}}
 \subfigure[$b_3=120$]{\includegraphics[width=0.47\textwidth,keepaspectratio]{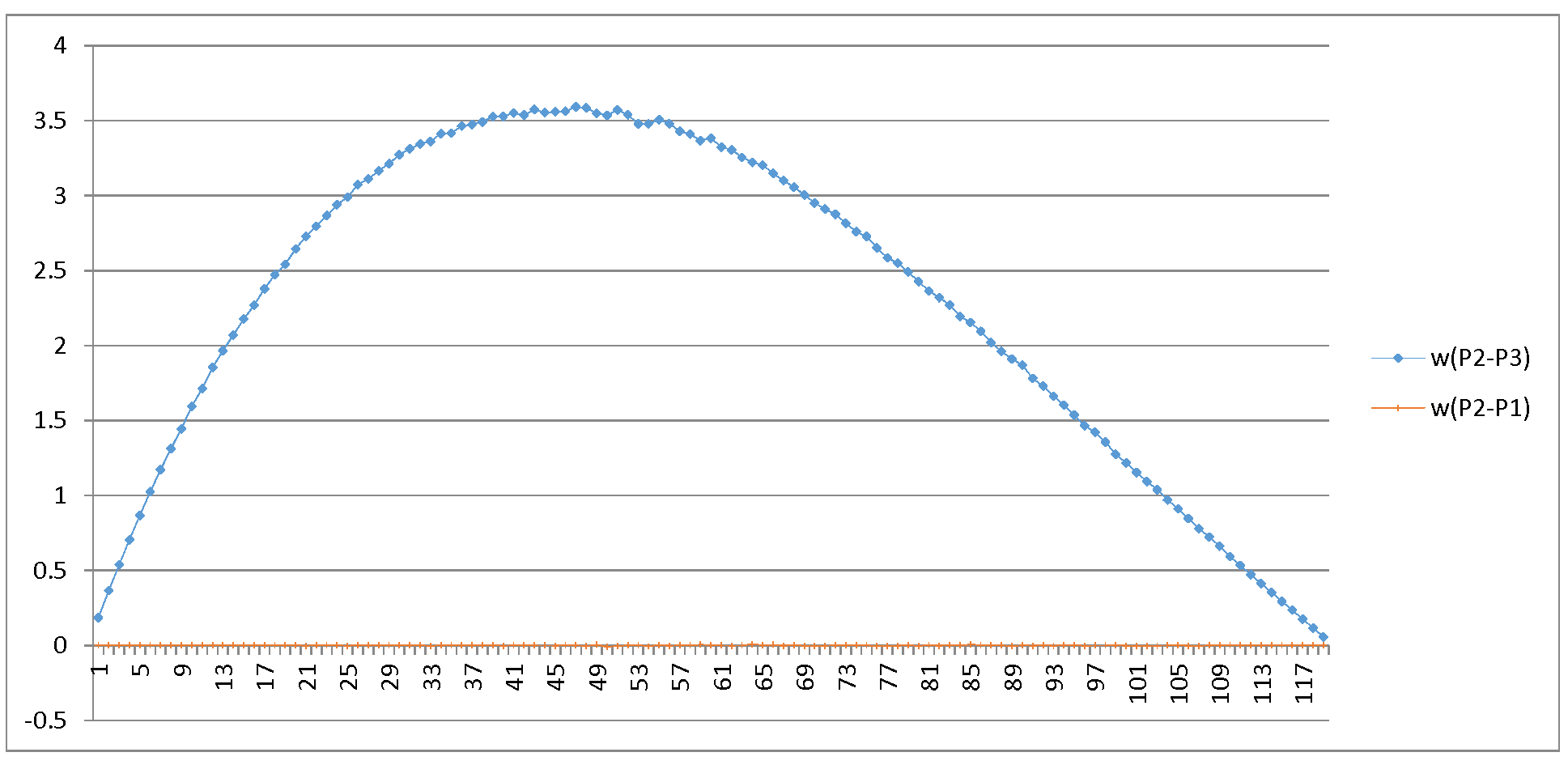}}
  \subfigure[$b_3=150$]{\includegraphics[width=0.47\textwidth,keepaspectratio]{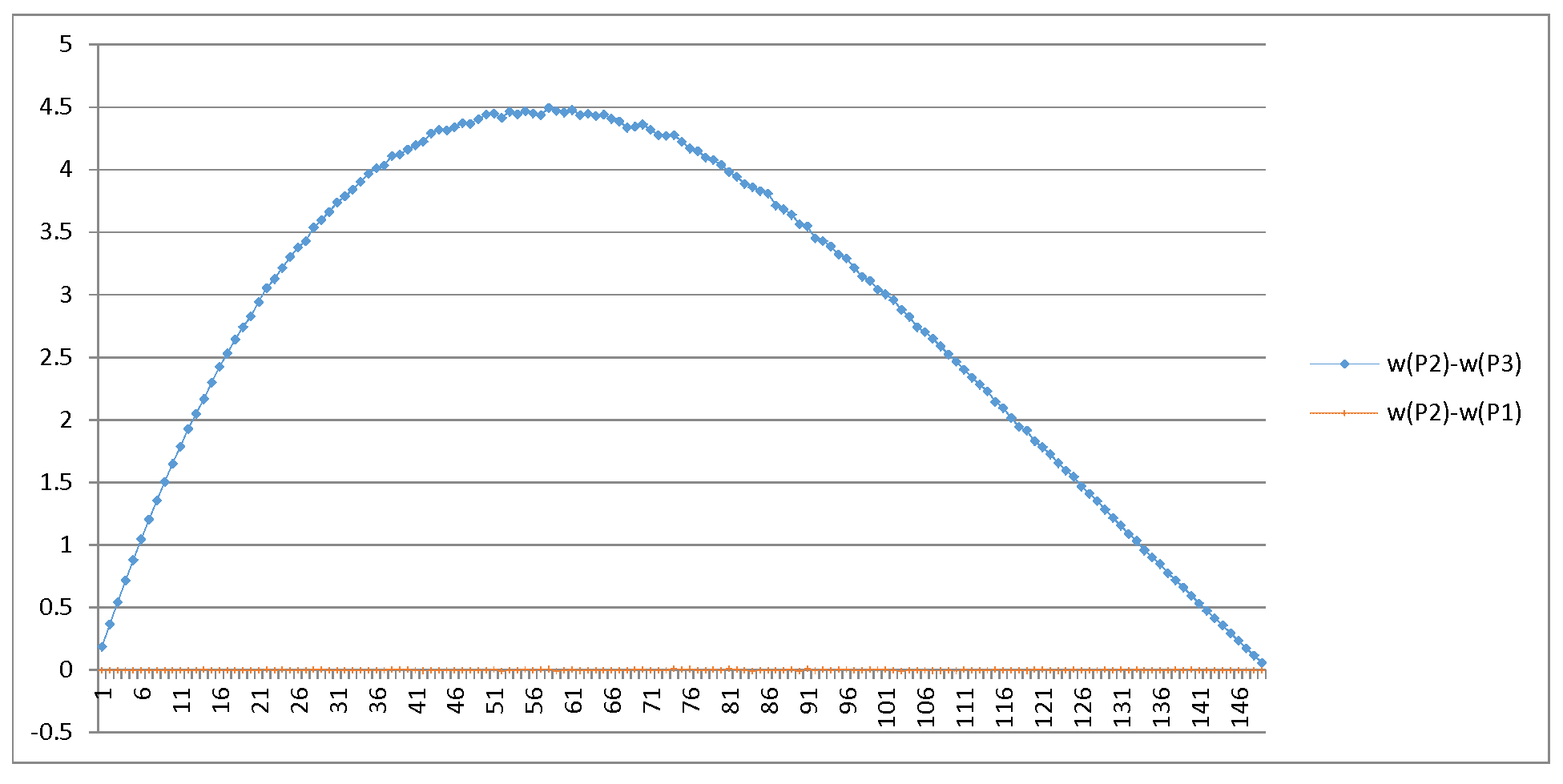}}
    \end{subfigmatrix}
    \caption{Worst-case analysis for $a_3$ ($a_1=a_2=0$, $b_1=b_2=1$)}
    \label{fig:worst}
 \end{figure}

Figure \ref{fig:worst} displays a plot of these results for $b_3=30,60,90,120 \; \text{and} \;150$.  From the inequality systems we know that $\Pol_h$ is exactly $\Pol_3$, therefore we are interested in the comparison between: $\Pol_2$ and $\Pol_3$, and $\Pol_2$ and $\Pol_1$.  From the plots of these differences, we see exactly what we would expect given the volume formulae.  The difference in mean width between $\Pol_2$ and $\Pol_1$ is very small; from a practical standpoint it is essentially zero.  The difference in mean width between $\Pol_2$ and $\Pol_3$ is always positive, indicating again that $\Po_3$ is the best choice of double-McCormick relaxation.  In addition, we observe that the maximum difference falls close to $a_3= \frac{b_3}{3}$ in all cases, demonstrating again that volume is a good predictor of how well a relaxation behaves.

%
%


%
%

%

\section*{Appendix}


\FloatBarrier

\clearpage

\bibliographystyle{alpha}
\bibliography{VolExperiments}
\end{document}